\documentclass[12pt]{article}

\usepackage[margin=1in]{geometry}
\usepackage{graphicx}
\usepackage{amsmath}
\usepackage{amsfonts}
\usepackage{amsthm}
\usepackage{authblk}
\usepackage{geometry,amssymb,graphicx,tensor,epigraph,amsmath,amsthm,xspace,listings,stmaryrd,float,tikz-cd,mathtools,afterpage}
\usepackage[T1]{fontenc}
\usepackage[utf8]{inputenc}
\usepackage[english]{babel}
\usepackage[disable]{todonotes}
\usepackage{hyperref}
\hypersetup{
    colorlinks=true,
    linkcolor=black,
    filecolor=magenta,
    urlcolor=blue,
    pdftitle={Complex Markets and Mean Field Games: Beyond Basic Models},
    pdfpagemode=FullScreen,
    }

\usepackage{bm}
\newcommand{\juany}[1]{}
\newcommand{\agus}[1]{}
\newcommand{\agusRta}[1]{}
\newcommand{\rafa}[1]{}

\newtheorem*{obs*}{Remark}

\newtheorem*{lema*}{Lemma}

\newtheorem*{remark*}{Remark}

\theoremstyle{definition}

\newcommand{\keywords}[1]{\textbf{\textit{Keywords---}} #1}

\begin{document}

\title{Complex Markets and Mean Field Games:\\Beyond Basic Models}

\author[1]{Agust\'in Mu\~noz Gonz\'alez}

\affil[1]{Departamento de Matem\'aticas, Facultad de Ciencias Exactas y Naturales, Universidad de Buenos Aires, Buenos Aires, Argentina}

\maketitle
\date{}

\begin{abstract}
This work builds on the theoretical frameworks presented in \cite{munoz2025liquiditypoolsmfg} and \cite{munoz2025liquiditypoolsmfgwcosts}, where the strategic interactions among traders in a constant-product market-making protocol were modelled using mean field games (MFG), first without transaction costs and then incorporating them. Here we present the formulation of a more complete model that integrates three types of agents: traders, liquidity providers (LPs), and arbitrageurs. While we do not establish existence results for this general model, the formulation identifies the main technical difficulties and lays the groundwork for future work. The LP acts as a dominating player in the sense of \cite{bensoussan2014mfgdominating}: its strategy influences the mean field distribution of the traders, and the equilibrium is sought as a solution to the coupled system of three problems that constitute a Major-Minor game. The arbitrageurs operate by solving the optimization problem of \cite{angeris2019analysis}, and their impact on the LP is captured through the loss-versus-rebalancing of \cite{milionis2022ammLVR}. The material in this article should be read as an open research proposal rather than a collection of closed results.
\end{abstract}

\hspace{10pt}
\clearpage
\tableofcontents
\clearpage

% Introduction
\section{Introduction}

In prior work \cite{munoz2025liquiditypoolsmfg,munoz2025liquiditypoolsmfgwcosts} we modelled the strategic interaction of traders in an AMM liquidity pool, first without transaction costs and then with fees, obtaining existence and characterization results for the equilibrium. All those models treat traders as the central actor, keeping the rest of the environment fixed. A natural limitation of that approach is that it ignores the strategic decisions of liquidity providers and the role of arbitrageurs, whose presence is structural in any real AMM.

In this article we present the formulation of a more complete model that integrates three types of agents---traders, liquidity providers (LPs), and arbitrageurs---in a constant-product liquidity pool with transaction costs. While we do not establish existence results for this general model, the formulation identifies the main technical difficulties and lays the groundwork for future work. The material in this article should be read as an open research proposal rather than a collection of closed results.

The first novel element of the model is the inclusion of arbitrageurs, who adjust the pool price with respect to the external price. Arbitrageurs act whenever they detect a gap between the price inside the pool and the price on an external market, forcing realignment and eliminating arbitrage opportunities. This has a direct impact on LPs through the phenomenon known as \textit{loss-versus-rebalancing} \cite{milionis2022ammLVR}, whereby the actions of arbitrageurs cause a potential loss for LPs, directly connecting these two types of agents in the system dynamics.

LPs, in turn, play a fundamental role by providing the liquidity required for pool operations. Unlike traders, LPs do not pay fees when adding or withdrawing liquidity, but their actions alter the pool invariant $k_t$, directly affecting the ETH and USDT reserves. In this extended model, LPs optimize over progressively measurable controls, although the path-dependence of the cost functional on the integrals of their past actions introduces a technical difficulty that requires a state augmentation (see Section~\ref{sec:juego_lp}). This optimization also affects traders, since it modifies the \textit{slippage} in transactions, thereby connecting LP actions to trader behavior.

Finally, traders optimize their inventories and move the pool price when they trade, affecting the ETH reserves and the pool invariant. Since the reserves incorporate the cumulative impact of both the LP and the arbitrageurs, traders must account for slippage in their decisions, completing the interaction cycle among the three agents of the system.

The article is organized as follows. In Section~\ref{sec:agentes} we present the three agents of the model with their individual dynamics. In Section~\ref{sec:interacciones} we describe the relationships among them: LP impact on price, slippage, loss-versus-rebalancing, and total pool reserves. In Section~\ref{sec:juego_completo} we formulate the LP game and the traders' game in the mean field limit. In Section~\ref{sec:mfg_major_minor} we instantiate the coupled Major-Minor system for this model. Sections~\ref{sec:dificultades} and \ref{sec:conclusion} discuss the technical difficulties and the proposed directions for future work.

% Agents
\section{Agents}
\label{sec:agentes}

In this section we present the three agents that interact in the model. We begin with the traders, revisiting the representation from \cite{munoz2025liquiditypoolsmfgwcosts}, and then introduce two new agents into the system, each with their own utility functions and control sets.

The first type of agent without whom liquidity pools would not exist is the liquidity provider (LP). This is the agent who holds both tokens available in the pool, ETH and USDT, in their inventory and is willing to deposit them into the pool in exchange for the fees that traders pay to swap tokens.

The second key agent is the arbitrageur. Their behavior is similar to that of a trader insofar as they operate in the pool, but with the difference that they possess information about prices on an external market and can detect and exploit the discrepancy between the pool price and that reference market \cite{angeris2019analysis}.

\subsection{Traders and pool}
Let $(\Omega,\mathcal{F}_T,\mathbb{F},\mathbb{P})$ be a filtered probability space supporting $N+1$ mutually uncorrelated Wiener processes $W_0,\dotsc,W_{N}$.

We revisit the trader dynamics, price dynamics, and pool reserve dynamics presented in \cite{munoz2025liquiditypoolsmfgwcosts}. We consider $N$ traders exchanging the two tokens ETH and USDT in a liquidity pool. Let $X_t^{i}$ and $Y_t^{i}$ denote the ETH and USDT inventories, respectively, of the $i$-th trader at time $t$. The dynamics of each inventory are given by

\begin{equation}
    \label{Pool 3: dinamica del trader}
    \begin{cases}
    \begin{aligned}
        dX_t^{i} &= \alpha_t^{i}\,dt + \sigma^{i}\,dW_t^i, \\
        dY_t^{i} &= -\alpha_t^i\,\frac{1+\phi^2}{2\phi}\,P_t\,dt,
    \end{aligned}
    \end{cases}
\end{equation}
where $\alpha_t^{i}:[0,T]\to \mathbb{R}$ represents the trading rate (the control); $\sigma^i$ is the volatility of trader $i$'s inventory, which for simplicity we assume to be independent of $i$; and $\sigma^i dW_t^i$ represents a random perturbation in the inventory. The USDT inventory dynamics are derived by assuming that the trader executes at the mid-price $\tilde{P}:=\frac{p^a+p^b}{2}$ (see \cite{mohan2021amm}), with $\phi = 1-\tau$ the fraction of assets received by the trader after paying the pool fee $\tau$.

If at time $t>0$ a trade of $\Delta_t^X$ units of ETH takes place, the price equation is written as

\begin{equation}
    \label{Pool 3: Eq: Ecuacion del precio base}
    P_t = \frac{k_0}{(X_0 + \phi\Delta_t^X)(X_0 + \Delta_t^X)},
\end{equation}
where $k_0 = X_0 Y_0$ is the initial pool invariant.

The dynamics of the ETH reserves inside the pool are modelled by the average of the traders' controls:
\begin{equation}
    \label{Pool 3: Eq: Cambio balance ETH base}
    X_t := X_0 - \frac{1}{N}\sum_{i=1}^{N}\int_0^t \alpha_s^{i}\,ds,
\end{equation}
where $X_0$ is the initial number of ETH tokens in the pool. We define the total inventory of the $i$-th trader at time $t$ as $V_t^i = Y_t^i + X_t^i P_t$.

We now introduce the working spaces used throughout the article.

\begin{itemize}
    \item Let $(\Omega,\mathcal{F},\mathbb{F}=(\mathcal{F}_t)_{0\leq t\leq T},\mathbb{P})$ be a complete filtered probability space, where the filtration $\mathbb{F}$ supports a one-dimensional Wiener process with respect to $\mathbb{F}$ and an initial condition $\xi\in L^2(\Omega,\mathcal{F}_0,\mathbb{P};\mathbb{R})$.

    \item Let $\mathcal{C} := C([0,T];\mathbb{R})$ be the space of continuous real-valued functions on $[0,T]$, equipped with the supremum norm $\|x\| := \sup_{s\in[0,T]}|x(s)|$.

    \item Given $\mathcal{P}(\mathbb{R})$ the space of probability measures on $\mathbb{R}$ and a measurable function $\psi:\mathbb{R}\to\mathbb{R}$, we define
    \begin{align*}
        \mathcal{P}_\psi(\mathbb{R}) &= \Bigl\{\mu\in\mathcal{P}(\mathbb{R}): \int\psi\,d\mu < \infty\Bigr\}, \\
        B_\psi(\mathbb{R}) &= \Bigl\{f:\Omega\to\mathbb{R}: \sup_\omega \frac{|f(\omega)|}{\psi(\omega)} < \infty\Bigr\}.
    \end{align*}
    We define $\tau_\psi(\mathbb{R})$ as the weakest topology on $\mathcal{P}_\psi(\mathbb{R})$ that makes the map $\mu\mapsto\int f\,d\mu$ continuous for every $f\in B_\psi(\mathbb{R})$.

    \item Let the control space $A\subset\mathbb{R}$ be a bounded subset, and let $\mathbb{A} = \{\alpha:[0,T]\times\Omega\to A: \text{progressively measurable}\}$ be the set of admissible controls. All controls are considered in closed-loop form.

    \item Let $\mathcal{P}(A)$ be the space of probability measures on $A$ together with the weak topology $\tau(A)$ that makes the map $q\mapsto\int_A f\,dq$ continuous for every $f\in B(A)$.
\end{itemize}

\subsection{Liquidity provider}
A liquidity provider deposits a certain amount of funds into the pool reserves. In exchange, they receive a form of receipt (typically in the form of share tokens) which they can later redeem for a fixed percentage of the reserves.

We consider a single LP as an agent who, at the start of the game, adds initial liquidity to the pool: $x_0$ units of ETH and $y_0$ units of USDT. Let $X_t^{LP}$ and $Y_t^{LP}$ denote the ETH and USDT inventories, respectively, held by the LP at time $t$. We write $Z_t^{LP}$ for the LP's pool share expressed in USDT given the liquidity provided.

\subsubsection{LP dynamics}

The dynamics of the ETH inventory are defined as
\begin{equation}\label{Pool 3: dinamica de X LP}
    dX_t^{LP} = \alpha_t^{LP}\,dt + \sigma^{X,LP}\,dW_t^{X,LP},
\end{equation}
where $\alpha_t^{LP}:[0,T]\to\mathbb{R}$ is the liquidity provision rate (the control); $\sigma^{X,LP}$ is the associated volatility; and $\sigma^{X,LP}dW_t^{X,LP}$ represents a random perturbation in the LP's inventory.

To understand the USDT inventory dynamics, we must account for the relationships between the quantities of each token that the LP adds to the pool. From the literature on liquidity provision in constant-product markets (see, e.g., \cite{angeris2019analysis}, \cite{heimbach2021lpdex}, \cite{zhang2018formal}), we know that the ratio between $\alpha_t^{LP}$ units of ETH and $\beta_t^{LP}$ units of USDT added by the LP must respect the current reserve ratio $P_t = Y_t/X_t$:
\begin{equation}
    \label{Pool 3: Eq: Relacion LP control ETH vs control USDT}
    \frac{\beta_t^{LP}}{\alpha_t^{LP}} = P_t \quad\Rightarrow\quad \beta_t^{LP} = \alpha_t^{LP} P_t.
\end{equation}
In other words, the LP must add liquidity without modifying the internal pool price. The USDT inventory dynamics are therefore
$$dY_t^{LP} = \beta_t^{LP}\,dt + \sigma^{Y,LP}\,dW_t^{Y,LP} = \alpha_t^{LP} P_t\,dt + \sigma^{Y,LP}\,dW_t^{Y,LP},$$
with $\sigma^{Y,LP}dW_t^{Y,LP}$ a random perturbation in the LP's USDT inventory.

Finally, at every point in time the LP's share of the pool is given by
\begin{equation}
    \label{Pool 3: Eq: Relacion LP control ETH vs control share}
    \gamma_t^{LP} := \frac{\alpha_t^{LP}}{X_t} = \frac{\beta_t^{LP}}{Y_t}.
\end{equation}
The second equality follows from \eqref{Pool 3: Eq: Relacion LP control ETH vs control USDT}.

For the controls to be consistent with the trader model and with the notion of change in the LP's inventory, we adopt the following convention:
\begin{itemize}
    \item $\alpha_t^{LP} > 0$ means the LP's ETH inventory increased, and hence their pool share decreases;
    \item $\alpha_t^{LP} < 0$ means the LP's ETH inventory decreased, and hence their pool share increases.
\end{itemize}
We allow the control to be negative. Under this convention, the relationship between $\alpha_t^{LP}$ and $\gamma_t^{LP}$ becomes
\begin{equation}
    \label{Pool 3: Eq: Relacion LP control ETH vs control share (Generalizada)}
    \gamma_t^{LP} = -\frac{\alpha_t^{LP}}{X_t}.
\end{equation}

The dynamics of the LP's pool share are given by
\begin{align*}
    dZ_t^{LP} = \gamma_t^{LP} V_t^P\,dt + \sigma^{Z,LP}\,dW_t^{Z,LP}
    = -\frac{\alpha_t^{LP}}{X_t} V_t^P\,dt + \sigma^{Z,LP}\,dW_t^{Z,LP},
\end{align*}
with $\sigma^{Z,LP}dW_t^{Z,LP}$ a random perturbation in the share and $V_t^P = 2Y_t = 2\sqrt{k_t P_t}$ the value of the pool in USDT units at time $t$ (see \cite{angeris2019analysis}, \cite{milionis2022ammLVR}).

In summary, the LP dynamics are given by
\begin{equation}
    \label{Pool 3: Dinamicas estados del LP}
    \begin{aligned}
        dX_t^{LP} &= \alpha_t^{LP}\,dt + \sigma^{X,LP}\,dW_t^{X,LP} && \longrightarrow \text{ETH holdings,} \\
        dY_t^{LP} &= \alpha_t^{LP} P_t\,dt + \sigma^{Y,LP}\,dW_t^{Y,LP} && \longrightarrow \text{USDT holdings,} \\
        dZ_t^{LP} &= -2\alpha_t^{LP} P_t\,dt + \sigma^{Z,LP}\,dW_t^{Z,LP} && \longrightarrow \text{pool share.}
    \end{aligned}
\end{equation}

Finally, we define
$$V_t^{LP} := X_t^{LP} P_t + Y_t^{LP} + Z_t^{LP}$$
as the total inventory of the LP at time $t$.

\subsection{Arbitrageur}
\subsubsection{Optimization problem}

We now turn to the third agent, the arbitrageurs. The arbitrageur's objective is to exploit the discrepancy between the market price $m_p$ and the price inside the pool $m_u$. To understand the full process and the no-arbitrage condition, we follow the logic developed in \cite{angeris2019analysis}.

Suppose the arbitrageur holds $\Delta_\beta$ USDT and wishes to exchange them for ETH. After applying the pool fee $\tau$, they exchange $\phi\Delta_\beta$ USDT for $\Delta_\alpha$ ETH at price $m_u$, i.e.,
$$\Delta_\alpha = \frac{\phi\Delta_\beta}{m_u},$$
with $\phi = 1-\tau$.

At some point they observe an opportunity on the external market and decide to sell those ETH at price $m_p$, receiving $\widetilde{\Delta}_\beta$ USDT with
$$\widetilde{\Delta}_\beta = \Delta_\alpha m_p = \phi\Delta_\beta\,\frac{m_p}{m_u}.$$
The trade is profitable whenever $\widetilde{\Delta}_\beta - \Delta_\beta > 0$, i.e., when

\begin{align*}
    \widetilde{\Delta}_\beta > \Delta_\beta
    &\Leftrightarrow \phi\,\frac{m_p}{m_u} > 1
    \Leftrightarrow (1-\tau)\,m_p > m_u.
\end{align*}

Therefore, the first no-arbitrage condition is $(1-\tau) m_p \leq m_u$. If one starts with ETH instead of USDT, one arrives at the second condition $m_u \leq (1+\tau)m_p$. The \textbf{no-arbitrage condition} is thus
\begin{equation}
    \label{Pool 3: Condicion de no arbitraje}
    (1-\tau)\,m_p \leq m_u \leq (1+\tau)\,m_p.
\end{equation}

The arbitrageur's problem can be written as the following optimization problem:
\begin{equation}
    \label{Pool 3: Arbitrator: Problema de optimizacion}
    \begin{cases}
    \begin{aligned}
    \text{maximize} \hspace{10pt} & m_p\Delta_\alpha - \Delta_\beta \\
    \text{subject to} \hspace{10pt} & \Delta_\alpha,\,\Delta_\beta \geq 0 \\
    & (R_\alpha - \Delta_\alpha)(R_\beta + \phi\Delta_\beta) = k.
    \end{aligned}
    \end{cases}
\end{equation}

Using the constant-product formula $X_t Y_t = k_t$ we can rewrite
$$\Delta_\beta = \frac{1}{\phi}\Bigl(\frac{k}{R_\alpha - \Delta_\alpha} - R_\beta\Bigr).$$

Since the map $x\mapsto 1/x$ is convex for $x>0$, $\Delta_\beta$ is a convex function of $\Delta_\alpha$. The problem equivalent to \eqref{Pool 3: Arbitrator: Problema de optimizacion} is
\begin{equation}
    \label{Pool 3: Arbitrator: Problema de optimizacion convexo}
    \begin{cases}
    \begin{aligned}
    \text{maximize} \hspace{10pt} & m_p\Delta_\alpha - \frac{1}{\phi}\Bigl(\frac{k}{R_\alpha - \Delta_\alpha} - R_\beta\Bigr) \\
    \text{subject to} \hspace{10pt} & \Delta_\alpha \geq 0.
    \end{aligned}
    \end{cases}
\end{equation}

In \cite{angeris2019analysis} it is shown that the optimal solution is the pair
\begin{equation}
    \label{Pool 3: solucion optima arbitrador}
    \begin{aligned}
    \Delta_\alpha^* &= R_\alpha - \sqrt{\frac{k}{\phi m_p}}, \\[4pt]
    \Delta_\beta^* &= \frac{1}{\phi}\Bigl(\sqrt{\phi m_p\, k} - R_\beta\Bigr),
    \end{aligned}
\end{equation}
provided $\phi m_p \leq m_u$. The arbitrageur solves this problem exogenously to the LP--trader game; their impact on the reserves and on the LP is incorporated in Section~\ref{sec:interacciones}.

% Relationships among agents
\section{Relationships among agents}
\label{sec:interacciones}

Now that we have described each agent in the system, their controls, and their dynamics, we establish the impact that each action has on the remaining participants and on the system as a whole. These interactions complicate the equations from the earlier models, and we therefore devote this section to developing them in detail.

\subsection{Impact of the LP on the price dynamics}

LPs provide liquidity without causing a change in price or paying fees. However, the equilibrium state of the pool is no longer given by the initial reserves $(X_0, Y_0)$; it must account for the cumulative impact of this agent. We define the new equilibrium states of the ETH and USDT reserves in the pool at time $t$ as
\begin{equation}
    \label{Eq: Reservas ajustadas por LP}
    \begin{aligned}
    X_t^{adj} &:= X_0 + \int_0^t \alpha_s^{LP}\,ds, \\
    Y_t^{adj} &:= Y_0 + \int_0^t \alpha_s^{LP} P_s\,ds.
    \end{aligned}
\end{equation}

Let us rewrite the price update equation taking into account this new equilibrium state $(X_t^{adj}, Y_t^{adj})$. In a liquidity pool with transaction costs, trades occur in two stages \cite{mohan2021amm}. Suppose that at time $t>0$ a trade of $\Delta_t^X$ units of ETH is generated; in the first stage the pool computes the number $\Delta_t^Y$ of USDT it must deliver in exchange, starting from the new equilibrium:
$$k_0 = (X_t^{adj} + \phi\Delta_t^X)(Y_t^{adj} - \Delta_t^Y) \quad\Rightarrow\quad \Delta_t^Y = Y_t^{adj} - \frac{k_0}{X_t^{adj} + \phi\Delta_t^X}.$$

In the second stage, fees are added to the pool and the invariant is updated:
$$k_t = X_t Y_t = (X_t^{adj} + \Delta_t^X)(Y_t^{adj} - \Delta_t^Y) = \frac{(X_t^{adj} + \Delta_t^X)\,k_0}{X_t^{adj} + \phi\Delta_t^X}.$$

The price equation then becomes
\begin{equation}
    \label{Pool 3: Eq: Nueva ecuacion del precio}
    P_t = \frac{k_0}{(X_t^{adj} + \phi\Delta_t^X)(X_t^{adj} + \Delta_t^X)}.
\end{equation}

Taking $\Delta_t^X = -\frac{1}{N}\sum_{i=1}^N\int_0^t\alpha_s^i\,ds$ as the average impact of traders, the differential dynamics of the price become
\begin{equation}
    \label{Pool 3: Eq: Dinamica del precio}
    \begin{aligned}
    dP_t &= \frac{\partial}{\partial t}\Biggl(\frac{k_0}{(X_t^{adj} + \phi\Delta_t^X)(X_t^{adj} + \Delta_t^X)}\Biggr)dt + \sigma_0\,dW_t^0 \\
    &= (-1)k_0\,
    \frac{2\Bigl(\alpha_t^{LP}X_t^{adj} + \phi\Delta_t^X\,\frac{1}{N}\sum_i\alpha_t^i\Bigr) + (1+\phi)\Bigl(X_t^{adj}\,\frac{1}{N}\sum_i\alpha_t^i + \alpha_t^{LP}\Delta_t^X\Bigr)}{\bigl((X_t^{adj} + \phi\Delta_t^X)(X_t^{adj} + \Delta_t^X)\bigr)^2}\,dt \\
    &\quad + \sigma_0\,dW_t^0,
    \end{aligned}
\end{equation}
where we have used $\frac{\partial}{\partial t}X_t^{adj} = \alpha_t^{LP}$ and $\frac{\partial}{\partial t}\Delta_t^X = \frac{1}{N}\sum_i\alpha_t^i$. A detailed derivation of this computation can be found in the appendix.

\subsection{Relationship between LP and traders: slippage}

The liquidity provider's action affects the traders' utility function through slippage. Slippage is the discrepancy between the market price and the effective execution price of a trade. In a constant-product AMM, the slippage experienced by a trader who exchanges $x_1$ units of token A for $x_2$ units of token B can be expressed as
$$S(x_1) = \frac{x_1/x_2}{m} - 1 = \frac{x_1}{r_1},$$
where $r_1$ is the balance of token 1 in the pool and $m$ is the spot price \cite{xu2022SoK}. This effect is amplified in smaller pools, where each trade significantly affects the relative quantities of assets.

To endow the model with greater realism and interconnection between the LP and traders, we propose incorporating $S$ as a slippage cost in the USDT inventory dynamics of trader $i$:
\begin{equation}\label{Pool 3: dinamica de X con slippage}
    \begin{aligned}
        dX_t^{i} &= \alpha_t^{i}\,dt + \sigma^i\,dW_t^i, \\
        dY_t^{i} &= -\alpha_t^i\,(1 - S(t,\alpha_t^i))\,\frac{1+\phi^2}{2\phi}\,P_t\,dt,
    \end{aligned}
\end{equation}
where $S(t,\alpha_t^i) = \frac{\alpha_t^i}{X_t}$. Note that slippage affects only what the trader receives in USDT, not what they deliver to the pool (which remains $\alpha_t^i$); the trader's ETH dynamics carry no slippage cost.

\subsection{Relationship between LP and arbitrageur: loss-versus-rebalancing}

There is a well-characterized relationship between the liquidity provider and the arbitrageur, described in \cite{milionis2022ammLVR}. In that article the authors decompose the value of the pool into a market-risk component and a non-negative, non-decreasing, predictable component they call \textit{loss-versus-rebalancing} (LVR).

This second component, obtained by replicating the movement of the pool balances in the external market as those balances change due to arbitrageur activity, can be interpreted as a cooperation cost that must be offset by fee revenue, an adverse-selection cost in the form of arbitrageur profits extracted from the pool, and an information cost because the pool lacks access to precise market prices.

Specifically, what the authors prove is that the pool value $V(P_t)$ can be written as
$$V(P_t) = R_t - LVR_t,$$
where, if the external market price follows a geometric Brownian motion that is a martingale under the risk-neutral measure
$$\frac{dP_t}{P_t} = \sigma\,dB_t^{\mathbb{Q}},$$
one defines:
\begin{itemize}
    \item The pool value function $V:\mathbb{R}_+\to\mathbb{R}_+$ given by the arbitrageur's optimization problem
    \begin{align*}
        V(P) := &\underset{(x,y)\in\mathbb{R}^2_+}{\text{minimize}} \hspace{6pt} Px + y, \\
        & \text{subject to} \hspace{6pt} f(x,y) = k;
    \end{align*}
    \item $\bm{R_t = V_0 + \int_0^t x^*(P_s)\,dP_s}$, the replication strategy for changes in pool balances on the external market, with $(x^*,y^*)$ the minimizers of the above problem;
    \item $\bm{LVR_t = \int_0^t \ell(P_s)\,ds}$, the loss-versus-rebalancing component, with
    $$\ell(P) = -\frac{\sigma^2 P^2}{2}\,V''(P) \geq 0.$$
\end{itemize}

Moreover, if the arbitrageur's gains are defined as $ARB_T = V(P_0) + \int_0^T x^*(P_t)\,dP_t - V(P_T)$, then $ARB_T = LVR_T$.

The computations in the cited article carry over analogously if one assumes active rather than passive LPs and transaction costs on trader actions. Our model is therefore consistent with the assumptions that derive the LVR formulation, and we can incorporate it into the LP's game as a running cost. We propose taking as the LP's total inventory function
\begin{equation}
    \label{Pool 3: Eq: Inventario LP - LVR}
    \widetilde{V}_t^{LP} = V_t^{LP} - LVR_t,
\end{equation}
whose dynamics are given by
\begin{equation}
    \label{Pool 3: Eq: dinamica Inventario LP - LVR}
    d\widetilde{V}_t^{LP} = dV_t^{LP} - \ell(P_t)\,dt.
\end{equation}

\subsection{Impact on pool reserves}

Following the logic of the models from earlier chapters, we model the ETH reserves inside the pool by the average of the traders' controls, with the difference that there is now an equilibrium state that depends on the LP's action and an impact from arbitrageurs who also pay fees when they trade. We define the ETH reserves inside the pool as
\begin{equation}
    \label{Pool 3: Eq: Reservas del pool con todas las interacciones}
    X_t := \underbrace{X_t^{adj}}_{\text{LP-induced equilibrium}} + \underbrace{\int_0^t \ell(P_s)\,ds}_{\text{arbitrageur impact}} - \underbrace{\frac{1}{N}\sum_{i=1}^N\int_0^t \alpha_s^i\,ds}_{\text{trader impact}}.
\end{equation}

% Full model formulation
\section{Model formulation}
\label{sec:juego_completo}

Having presented all the final relationships and dynamics, we now state the game solved by the LP and the game solved by the traders. Based on the impact on pool reserves analyzed in the previous section, we write
$$\Delta_t^X = \int_0^t \ell(P_s)\,ds - \frac{1}{N}\sum_{i=1}^N\int_0^t\alpha_s^i\,ds.$$

We also define
$$H^q(t) := \int_0^t\ell(P_s)\,ds - \int_0^t\int_A \tilde{a}\,dq_s(\tilde{a})$$
as the total impact on reserves in the limit as the number of traders tends to infinity, where $q_t$ is the law of the traders' controls at time $t$.

\subsection{LP game}
\label{sec:juego_lp}
To write the LP's optimization problem, we compute the dynamics of their inventory given the price equation \eqref{Pool 3: Eq: Nueva ecuacion del precio} and the reserves \eqref{Pool 3: Eq: Reservas del pool con todas las interacciones}. Applying It\^{o}'s lemma to the adjusted inventory \eqref{Pool 3: Eq: Inventario LP - LVR} we obtain

\begin{equation}
    \label{Pool 3: Ito_inventario_LP}
    \begin{aligned}
    d(V_t^{LP})
        &= d(X_t^{LP} P_t) + d(Y_t^{LP}) + d(Z_t^{LP}) \\
        &= (-1)k_0\,X_t^{LP}
        \frac{2\Bigl(\alpha_t^{LP}X_t^{adj} + \phi\Delta_t^X\,\frac{1}{N}\sum_i\alpha_t^i\Bigr) + (1+\phi)\Bigl(X_t^{adj}\,\frac{1}{N}\sum_i\alpha_t^i + \alpha_t^{LP}\Delta_t^X\Bigr)}{\bigl((X_t^{adj} + \phi\Delta_t^X)(X_t^{adj} + \Delta_t^X)\bigr)^2}\,dt \\
        &\quad + \sigma^{X,LP} P_t\,dW_t^{X,LP} + X_t^{LP}\sigma_0\,dW_t^0 + \sigma^{Y,LP}\,dW_t^{Y,LP} + \sigma^{Z,LP}\,dW_t^{Z,LP}.
    \end{aligned}
\end{equation}

We assume the LP is risk-neutral and seeks to maximize their expected profit from operating in the decentralized market:
$$J^{LP}(\alpha^{LP}) = \mathbb{E}\Biggl[\int_0^T f^{LP}(t,X_t,\alpha_t)\,dt - l(X_T,Z_T)\Biggr],$$
where $l:\mathbb{R}\to\mathbb{R}$ represents a terminal inventory cost. By equation \eqref{Pool 3: Ito_inventario_LP} and the LVR dynamics \eqref{Pool 3: Eq: dinamica Inventario LP - LVR}, the function $f^{LP}:[0,T]\times\mathcal{C}\times\widetilde{A}\to\mathbb{R}$ is defined by

\begin{align*}
    f^{LP}(t,x,\alpha)
    &= -k_0\,x_t\,G(t,\alpha)^2\Bigl[2\Bigl(\alpha_t\,x_{\alpha,t}^{adj} + \phi\,H^q(t)\textstyle\int_A\tilde{a}\,dq_t(\tilde{a})\Bigr) \\
    &\quad + (1+\phi)\Bigl(x_{\alpha,t}^{adj}\textstyle\int_A\tilde{a}\,dq_t(\tilde{a}) + \alpha_t\,H^q(t)\Bigr)\Bigr],
\end{align*}
where $q$ is the law of the traders' controls; $x_{\alpha,t}^{adj} := X_0 + \int_0^t\alpha_s\,ds$; and
$$G(t,\alpha) := \frac{1}{(x_{\alpha,t}^{adj} + H^q(t))(x_{\alpha,t}^{adj} + \phi\,H^q(t))}.$$

\begin{obs*}[Path-dependence and LP control space]
The function $f^{LP}$ involves the cumulative quantity $x_{\alpha,t}^{adj} = X_0 + \int_0^t\alpha_s^{LP}\,ds$. To maintain a Markovian formulation, it suffices to introduce the auxiliary state $S_t := \int_0^t\alpha_s^{LP}\,ds$, so that $X_t^{adj} = X_0 + S_t$ and $dS_t = \alpha_t^{LP}\,dt$. With this augmented state, the LP optimizes over the standard space of progressively measurable controls without the need to impose continuity on the control paths.
\end{obs*}

\subsection{Traders' game}
\label{sec:juego_traders}
Taking into account the interactions presented in the previous section, we now derive the trader's problem. By It\^{o}'s lemma we have
\begin{equation}
\label{Pool 3: Ito_inventario_trader}
    \begin{aligned}
    dV_t^i &= dY_t^i + X_t^i\,dP_t + dX_t^i\,P_t \\
    &= \Biggl[X_t^i\,(-1)k_0
    \frac{2\Bigl(\alpha_t^{LP}X_t^{adj} + \phi\Delta_t^X\,\frac{1}{N}\sum_i\alpha_t^i\Bigr) + (1+\phi)\Bigl(X_t^{adj}\,\frac{1}{N}\sum_i\alpha_t^i + \alpha_t^{LP}\Delta_t^X\Bigr)}{\bigl((X_t^{adj} + \phi\Delta_t^X)(X_t^{adj} + \Delta_t^X)\bigr)^2} \\
    &\quad + \alpha_t^i\,\frac{k_0}{(X_t^{adj} + \phi\Delta_t^X)(X_t^{adj} + \Delta_t^X)} \\
    &\quad + \alpha_t^i\,(1 - S(t,\alpha_t^i))\Bigl(1 - \frac{1+\phi^2}{2\phi}\Bigr)\frac{k_0}{(X_t^{adj} + \Delta_t^X)(X_t^{adj} + \phi\Delta_t^X)}\Biggr]dt \\
    &\quad + X_t^i\sigma_0\,dW_t^0 + P_t\sigma^i\,dW_t^i.
    \end{aligned}
\end{equation}

Assuming agents are risk-neutral and seek to maximize their expected profit, the $i$-th trader seeks to maximize
$$J^i(\alpha^1,\dotsc,\alpha^N) = \mathbb{E}\Biggl[\int_0^T f(t,X_t^i,\hat\nu_t^N,\alpha_t^i)\,dt - l(X_T^i)\Biggr],$$
where $\hat\nu_t^N$ denotes the empirical distribution of $\alpha_t^1,\dotsc,\alpha_t^N$, $l:\mathbb{R}\to\mathbb{R}$ represents a terminal inventory cost, and, by equation \eqref{Pool 3: Ito_inventario_trader}, the function $f:[0,T]\times\mathcal{C}\times\mathcal{P}_\psi(\mathcal{C})\times\mathcal{P}(A)\times A\to\mathbb{R}$ is defined by

\begin{align*}
    f(t,x,\mu,q,\alpha)
    &= -k_0\,x_t\,G(t,\alpha^{LP})^2\Bigl[2\Bigl(\alpha_t^{LP}X_t^{adj} + \phi\,H^q(t)\textstyle\int_A\tilde{a}\,dq_t(\tilde{a})\Bigr) \\
    &\quad + (1+\phi)\Bigl(X_t^{adj}\textstyle\int_A\tilde{a}\,dq_t(\tilde{a}) + \alpha_t^{LP}\,H^q(t)\Bigr)\Bigr] \\
    &\quad + \alpha\,k_0\,G(t,\alpha^{LP}) + \alpha\,k_0\,G(t,\alpha^{LP})\,(1 - S(t,\alpha))\Bigl(1 - \frac{1+\phi^2}{2\phi}\Bigr),
\end{align*}
for $X_t^{adj}$ and $G(t,\alpha^{LP})$ defined from the history of the controls $\alpha^{LP}$.

By the symmetry of the model, trader $i$'s contribution to $\hat\nu_N$ is negligible for large $N$, and we may treat $\mu,q$ as fixed. The objective in the mean field limit is therefore
$$J(\alpha) := \mathbb{E}\Biggl[\int_0^T f(t,X_t,\mu_t,q_t,\alpha_t)\,dt - l(X_T)\Biggr],$$
where $(\mu_t)_{0\leq t\leq T}$ is a flow of probability measures on the state space and $(q_t)_{0\leq t\leq T}$ is a flow of measures on the control space.

% Major-Minor formulation
\section{Major-Minor formulation}
\label{sec:mfg_major_minor}

Following the framework of \cite{bensoussan2014mfgdominating}, we consider a mean field game between the dominating player (LP) and a group of representative agents (traders), each of whom behaves similarly and also interacts with the others through a mean field term substantially influenced by the dominating player. Unlike a classical Stackelberg game in which the leader commits sequentially, the structure here is of Major-Minor type: the LP influences the mean field distribution of the traders, and the equilibrium is sought simultaneously for the three coupled problems.

The representative agent's control problem must first be solved, after which the equilibrium condition leads to coupled Hamilton-Jacobi-Bellman and Fokker-Planck equations. Because the mean field term is endogenous to the dominating player, in order to achieve optimal control the latter must account for the coupled equations when deciding their own strategy.

We use the notation of the General Major-Minor framework and write the mean field term through the joint state--control law
\[
\Pi_t = \mathcal{L}\bigl(X_t^1, U_t^1 \,\big|\, \mathcal{F}_t^0\bigr)\in\mathcal{P}_2(\mathbb{R}^{n_1}\times A),
\]
which describes the conditional distribution of the representative agent's pair $(x_1,u_1)$ given the common noise of the dominating player. The law of states and the law of controls are recovered as marginals:
\[
\mu_t = (\mathrm{pr}_X)_{\#}\Pi_t, \qquad q_t = (\mathrm{pr}_U)_{\#}\Pi_t.
\]
We use the notation $x_0,u_0$ for the state and control of the dominating player (LP), and $x_1,u_1$ for the state and control of the representative agent (trader).

\subsection{Problem 1: Representative agent control}
\label{Problema final 1}

Given the process $x_0$ and an exogenous probability measure $\nu$ (adapted to $\mathcal{F}_t^0$), find a control $u_1\in\mathcal{A}_1$ that minimizes the cost functional
\[
J_1(x, \alpha, q) := \mathbb{E}\Biggl[\int_0^T f_1\bigl(x_1(t), x_0(t), \nu(t), u_1(t)\bigr)\,dt + h_1\bigl(x_1(T), x_0(T), \nu(T)\bigr)\Biggr],
\]
where the function $f_1:[0,T]\times\mathcal{C}\times\mathcal{P}_\psi(\mathcal{C})\times\mathcal{P}(A)\times A\to\mathbb{R}$ is defined by
\begin{align*}
    f_1(t,x,\mu,q,\alpha)
    &= -k_0\,x_1(t)\,G(t,\alpha^{LP})^2\Bigl[2\Bigl(\alpha_t^{LP}X_t^{adj} + \phi\,H^q(t)\textstyle\int_A\tilde{a}\,dq_t(\tilde{a})\Bigr) \\
    &\quad + (1+\phi)\Bigl(X_t^{adj}\textstyle\int_A\tilde{a}\,dq_t(\tilde{a}) + \alpha_t^{LP}\,H^q(t)\Bigr)\Bigr] \\
    &\quad + \alpha\,k_0\,G(t,\alpha^{LP}) + \alpha\,k_0\,G(t,\alpha^{LP})\,(1 - S(t,\alpha))\Bigl(1 - \frac{1+\phi^2}{2\phi}\Bigr),
\end{align*}
and we take $h_1(x_1(T),x_0(T),\nu(T)) = l(x_1(T))$.

\subsection{Problem 2: Mean field consistency condition}
\label{Problema final 2}

Given an exogenous process-valued probability measure $\nu$, let $\mathcal{M}(\nu)(t)$ be the measure induced by the optimal state $x_1(t)$ found in Problem~\ref{Problema final 1} conditioned on $\mathcal{F}_t^0$. Find the process-valued probability measure $\mu$ satisfying the fixed-point property: $\mathcal{M}(\mu)(\cdot) = \mu(\cdot)$.

\subsection{Problem 3: Dominating player control}
\label{Problema final 3}

Find a control $u_0\in\mathcal{A}_0$ that minimizes the cost functional
\[
J_0(u_0) := \mathbb{E}\Biggl[\int_0^T f_0\bigl(x_0(t), \mu(t), u_0(t)\bigr)\,dt + h_0\bigl(x_0(T), \mu(T)\bigr)\Biggr],
\]
where $\mu$ is the solution found in Problem~\ref{Problema final 2} and the function $f^{LP}:[0,T]\times\mathcal{C}\times\widetilde{A}\to\mathbb{R}$ is defined by

\begin{align*}
    f^{LP}(t,x,\alpha)
    &= -k_0\,x_t\,G(t,\alpha)^2\Bigl[2\Bigl(\alpha_t\,x_{\alpha,t}^{adj} + \phi\,H^q(t)\textstyle\int_A\tilde{a}\,dq_t(\tilde{a})\Bigr) \\
    &\quad + (1+\phi)\Bigl(x_{\alpha,t}^{adj}\textstyle\int_A\tilde{a}\,dq_t(\tilde{a}) + \alpha_t\,H^q(t)\Bigr)\Bigr].
\end{align*}

Note that $f^{LP}$ depends on the joint law $\Pi_t$ through its marginals: the law of states $\mu_t$ (implicitly through the pool reserves) and the law of controls $q_t$ (explicitly through the integral $\int_A\tilde{a}\,dq_t(\tilde{a})$).

Problems~1--3 constitute the instantiation of the coupled SHJB--Fokker-Planck--Adjoint system of the General Major-Minor framework for the particular case of liquidity pools with three agents. The necessary optimality conditions and the complete structure of the coupled equation system are those of the main theorem of that framework, with the cost functions $f_1$ and $f^{LP}$ defined above.

% Technical difficulties and proposed directions
\section{Technical difficulties}
\label{sec:dificultades}

Establishing existence results for the full model faces several difficulties that are not present in the models of earlier chapters.

\begin{enumerate}
    \item \textbf{Multiplicative nonlinear coupling.} The traders' reward function depends on the product of their state, aggregate actions, pool reserves, and the liquidity provided by LPs. This multiplicative coupling prevents the additive decomposition required by the standard sandwich-type conditions and complicates the direct application of the existence techniques developed for the pool model with transaction costs.

    \item \textbf{Path-dependence.} The pool reserves depend on the complete history of all agents' actions (through the integrals $\int_0^t\alpha_s\,ds$ and the cumulative impact $H^q(t)$), introducing a trajectory dependence that goes beyond the standard Markovian setting.

    \item \textbf{Major-Minor game with mean field.} The combination of a dominating player (LP) with a mean field game among traders requires the simultaneous resolution of the three coupled problems of the SHJB--Fokker-Planck--Adjoint system: the representative agent's optimal control (Problem~\ref{Problema final 1}), the mean field consistency condition (Problem~\ref{Problema final 2}), and the dominating player's optimal control (Problem~\ref{Problema final 3}). The existence of equilibria for this coupled system, with the specific nonlinearities of AMM price functions, is an open problem.

    \item \textbf{Regularity of the Lions derivative.} The dependence on the joint state--control law $\Pi_t$ of the Major-Minor framework requires verifying that the Lions derivative $\frac{\partial f}{\partial\mu}$ is well defined under the specific nonlinearities of AMM price functions.
\end{enumerate}

\section{Proposed directions}
\label{sec:direcciones}

\paragraph{Tractable simplifications.} A first step is to solve the model under simplifications that reduce the technical difficulties: exogenous arbitrageur (external price as a given process), LP with fixed strategy (constant liquidity), or traders with a short horizon (myopic approximation). These reduced versions may be accessible using the techniques developed in \cite{munoz2025liquiditypoolsmfg,munoz2025liquiditypoolsmfgwcosts}.

\paragraph{Numerical validation.} The numerical pipeline developed for the base model (fixed-point solver, $\varepsilon$-Nash test, convergence in $N$) can be adapted to the full model to explore numerically the existence and stability of equilibria before seeking formal proofs.

\paragraph{Extension of the LVR analysis.} In this work we incorporate the results of \cite{milionis2022ammLVR} on loss-versus-rebalancing as a component of the model, treating arbitrageurs as agents who extract LVR from LPs. A natural direction is to analyze how the equilibrium of the Major-Minor game quantitatively modifies the LVR relative to the passive-LP case, and whether the LP's strategic optimization can reduce arbitrage losses.

\paragraph{Calibration with on-chain data.} Calibrating the model with real transaction data from Uniswap or other protocols would allow one to validate the qualitative predictions and to quantify the equilibrium effects in real markets.

% Conclusion
\section{Conclusion}
\label{sec:conclusion}

In this article we have presented the formulation of an AMM liquidity pool model with three types of agents: traders, liquidity providers, and arbitrageurs. The game structure is of Major-Minor type in the sense of \cite{bensoussan2014mfgdominating}: the LP acts as a dominating player whose strategy influences the mean field distribution of the traders, while the arbitrageurs operate exogenously by solving the problem of \cite{angeris2019analysis}, with their impact on the LP captured through the LVR of \cite{milionis2022ammLVR}.

The formulation identifies with precision the sources of model complexity: the multiplicative coupling in the reward functions, the path-dependence of the LP's controls on the pool reserves, and the need to simultaneously solve three coupled problems (representative agent control, mean field consistency, and dominating player control) under the nonlinearities of the AMM price function.

As noted in Section~\ref{sec:dificultades}, the existence of equilibria for this coupled system is an open problem. This article should be read as a research proposal that identifies the correct framework and the concrete technical difficulties, laying the groundwork for future work aimed both at analytical results under tractable simplifications and at numerical validation of the full model.

% Appendix
% \appendix
% \input{Parts/Appendix}

\clearpage
\bibliographystyle{plain}
\bibliography{main}

\end{document}